\documentstyle[12pt]{article}
\begin{document}

\begin{center}
{\Large \bf  Uncertainty Relations and Wave Packets on the Quantum Plane }
\end{center}

\vspace{1cm}

\begin{flushleft}
H. Ahmedov$^1$  and I. H. Duru$^{2,1}$

\vspace{.5cm}

{\small
1. TUBITAK, Feza G\"ursey Institute,  P.O. Box 6, 81220 Cengelkoy, 
Istanbul, Turkey 
\footnote{E--mail : hagi@fge1.gursey.gov.tr and duru@fge1.gursey.gov.tr}.

2. Trakya University, Mathematics Department, P.O. Box 126, 
Edirne, Turkey.}
\end{flushleft}

\vspace{2cm}

\begin{center}
{\bf \small  Abstract}
\end{center} 
(2+2)--dimensional quantum mechanical q--phase space which is the
semi--direct product of the quantum plane $E_q(2)/U(1)$ and its dual algebra 
$e_q(2)/u(1)$ is constructed. Commutation and the resulting uncertainty
relations are studied. "Quantum mechanical q--Hamiltonian'' of the motion
over the quantum plane is derived and the solution of the Schr\"{o}dinger
equation for the q--semiclassical motion governed by the expectation value
of that Hamiltonian is solved.

\vspace{2cm}

\begin{center}
July 1998
\end{center}

\pagebreak

\vfill
\eject

\renewcommand{\theequation}{I.\arabic{equation}} \setcounter{equation}{0}

\section{Introduction}

How the quantum mechanical effects are altered if we replace the space--time
continuum by a non--commutative geometry is an interesting question.
Especially the phenomena which are directly related to the geometry of the
environment, such as the Casimir effect, Aharonov--Bohm effect and the
cosmological pair production may exhibit peculiarities. To obtain the
required mathematical tools to study the above mentioned physical effect a
program dealing with the Green functions for q--group spaces has already
started \cite{1}, \cite{2}.

Another consequence of the non--commutative nature of the geometry is the
changes occurring in the structure of the quantum mechanical phase space:
What are the new commutation relations and the uncertainty relations? There
exists an extensive literature on the commutation relations between the
position and momentum operators resulting from the quantum group related
associative algebras and the corresponding new uncertainty relations \cite{3}%
. Since most of the above mentioned literatures employ 1--dimensional
position operators the corresponding configuration spaces are commutative.
The operators act on the Hilbert space of the square integrable functions of
the commutative coordinates or momenta.

In this paper we study the case of a 2--dimensional non--commutative
configuration space $E_q^2$ which is obtained from the space of q--Euclidean
group $E_q(2)$. We have mutually non--commutative position operators $\hat {%
{\bf x}}$ and $\hat {{\bf y}}$ ( and thus mutually non--commutative momentum
operators $\hat {{\bf p}}_x$, $\hat {{\bf p}}_y$). They act on the Hilbert
space $H(E_q^2)$ of square integrable functions of non--commutative
coordinates ${\bf x}$ and ${\bf y}$. The latter can be realized in the
Hilbert space $H(S)$ of square integrable functions on the circle $S$.

In the following section we present our ( 2+ 2)--dimensional quantum
mechanical q--phase space.

In Section III the wave packets are obtained and the uncertainty relations
are studied.

In Section IV the free motion Hamiltonian $\hat{{\bf H}}$ which commutes
with $\hat{{\bf p}}_{x}$, $\hat{{\bf p}}_y$ and the angular momentum like
operator $\hat{{\bf l}}$ is introduced. We then write the q--expectation
value of this Hamiltonian, and then considering the expectation value of $%
\hat{p}_x$, $\hat{p}_y$ and $\hat{{\bf l}}$ as continuously varying numbers
obtain a q--semiclassical Hamiltonian $H$. Schr\"{o}dinger equation
corresponding to this q--semiclassical Hamiltonian is also solved.

All the required background on the quantum groups and algebras are presented
in the Appendix.

\renewcommand{\theequation}{II.\arabic{equation}} \setcounter{equation}{0}

\section{(2+2)--Dimensional Quantum Mechanical q--Phase Space}

We start with defining the basic concepts that we employ. All the known
definitions related to the quantum group $E_q(2)$ and it's dual $e_q(2)$ are
given in the Appendix.

\vspace{1cm}

(i) q--configuration space $E_q^2$ is the subspace of the quantum group $%
E_q(2)$ defined as 
\begin{equation}
E_q^2=\sum_{j\in Z}\oplus E[j,0],
\end{equation}
where $E[i,j]$ is given by (\ref{eq: a7}). In other words q--configuration
space is the left sided coset space $E_q(2)/U(1)$. One can show that $E_q^2$
is generated by ${\bf z}_{\pm }$, which satisfy the commutation relation 
\begin{equation}
\label{eq: II0}[{\bf z}_{+},{\bf z}_{-}]=(1-q^2){\bf z}_{+}{\bf z}_{-}.
\end{equation}
We can also define the cartesian ${\bf x}$, ${\bf y}$ and the polar ${\bf r}$%
, ${\bf o}$ coordinates as 
\begin{equation}
\label{eq: II1}{\bf x}=\frac 1{\sqrt{2}}({\bf z}_{+}+{\bf z}_{-}),\ \ \ \ \
\ {\bf y}=\frac 1{\sqrt{2}i}({\bf z}_{+}-{\bf z}_{-})
\end{equation}
and 
\begin{equation}
\label{eq: II2}{\bf z}_{+}={\bf r}e^{i{\large {\bf o}}},\ \ \ \ {\bf z}%
_{-}=e^{-i{\bf o}}{\bf r}
\end{equation}
respectively. By the virtue of (\ref{eq: II0})  we get the
commutation relations 
\begin{equation}
\label{eq: II3}[{\bf x},{\bf y}]=i\tanh \Lambda ({\bf x}^2+{\bf y}^2),
\end{equation}
\begin{equation}
\label{eq: II4}[{\bf r},{\bf o}]=i\Lambda {\bf r},
\end{equation}
where $\Lambda =\log q^{-1},\ \ q<1$. Involutions are 
\begin{equation}
\label{eq: II5}{\bf x}^{*}={\bf x},\ \ \ \ {\bf y}^{*}={\bf y}.
\end{equation}
\begin{equation}
\label{eq: II6}{\bf r}^{*}={\bf r},\ \ \ \ {\bf o}^{*}={\bf o}.
\end{equation}
Inspecting the commutation relations (\ref{eq: II3}) and (\ref{eq: II4}) we
conclude that for larger radial distance one observes larger
non--commutativity. The operators ${\bf z}_{\pm }$ and consequently ${\bf x}$%
, ${\bf y}$ and ${\bf r}$, ${\bf o}$ are defined in the Hilbert space $H(S)$
(\ref{eq: b0}). In a fashion parallel to the quantum mechanics we can
interpret ${\bf x}$, ${\bf y}$ as observables on the states in $H(S)$. The
commutative coordinates $x$ and $y$ are the expectation values of the
self--adjoint operators ${\bf x}$ and ${\bf y}$ : 
\begin{equation}
x=(u,{\bf x}u),\ \ \ y=(u,{\bf y}u),
\end{equation}
where $u\in H(S)$.

\vspace{1cm}

(ii) q--momentum space $e_q^2$ is the quotient algebra $e_q(2)/u(1)$. It is
generated by ${\bf p}_\pm$, which by the virtue of (\ref{eq: a9}) satisfy
the commutation relation 
\begin{equation}
\label{eq: II7}[{\bf p}_+, {\bf p}_-]=(1-q^2){\bf p}_+ {\bf p}_- . 
\end{equation}
q--momentum space $e_q^2$ and q--configuration space $E_q^2$ are in duality.
We define the momenta ${\bf p}_x$, ${\bf p}_y$ in cartesian coordinates as 
\begin{equation}
\label{eq: II8}{\bf p}_x=\frac{1}{\sqrt{2}} ({\bf p}_+ +{\bf p}_- ), \ \ \ \
\ \ {\bf p}_y=\frac{i}{\sqrt{2}} ({\bf p}_+ -{\bf p}_- ). 
\end{equation}
By the virtue of (\ref{eq: II7}) and (\ref{eq: a10}) we get the commutation
relation 
\begin{equation}
\label{eq: II9}[{\bf p}_x, {\bf p}_y] =i \tanh \Lambda ({\bf p}_x^2 + {\bf p}%
_y^2 ) 
\end{equation}
and the involution 
\begin{equation}
{\bf p}_x^* ={\bf p}_x, \ \ \ \ {\bf p}_y^*={\bf p}_y. 
\end{equation}
We observe that ${\bf p}_x$ and ${\bf p}_y$ satisfy the same commutation
relations as ${\bf x}$ and ${\bf y}$. Therefore the elements of the
q--momentum space are operators in the Hilbert space $H(S)$. In the
orthonormal basis (\ref{eq: b00}) we have 
\begin{equation}
\label{eq: II}{\bf p}_+ e_j = p_0 q^{-j+\frac{1}{2}} e_{j-1}, \ \ \ {\bf p}%
_- e_j = p_0 q^{- j-\frac{1}{2}} e_{j+1}, 
\end{equation}
where $p_0$ is a constant of dimension of momentum. The commutative momenta $%
p_x$, $p_y$ are the expectation values of the operators ${\bf p}_x$ and $%
{\bf p}_y$ in $H(S)$.

q--phase space $E_q^2\otimes e_q^2$ is the direct product of
q--configuration and q--momentum spaces. It is realized in the Hilbert space 
$H(S)\otimes H(S)$.

\vspace{1cm}

(iii) ''quantum mechanical q--phase space'' $E_q^2\star e_q^2$ is the
semi--direct product of the q--configuration and q--momentum spaces. The
position $\hat {{\bf z}}_{\pm }$ and momentum $\hat {{\bf p}}_{\pm }$
operators are defined in the Hilbert space $H(E_q^2)$ as 
\begin{equation}
\label{eq: II10}\hat {{\bf z}}_{\pm }f=f{\bf z}_{\pm },\ \ \ \hat {{\bf p}%
}_{\pm }f=\hbar {\cal R}({\bf p}_{\pm })f,
\end{equation}
where $f\in \Phi (E_q^2)\subset H(E_q^2)$ and ${\cal R}$ is the right
representation (\ref{eq: b7}) of the quantum algebra $e_q(2)$ corresponding
to the left quasi--regular representation (\ref{eq: b6}) of the quantum
group $E_q(2)$. Here $\Phi (E_q^2)$ is dense subset in $H(E_q^2)$ defined in
(\ref{eq: b1}). The position $\hat {{\bf x}}$, $\hat {{\bf y}}$ and momentum 
$\hat {{\bf p}}_x$, $\hat {{\bf p}}_y$ operators in cartesian coordinates
are given by 
\begin{equation}
\label{eq: II11}\hat {{\bf x}}f=f{\bf x},\ \ \ \hat {{\bf y}}f=f{\bf y},
\end{equation}
\begin{equation}
\label{eq: II12}\hat {{\bf p}}_xf=\hbar {\cal R}({\bf p}_x)f,\ \ \ \hat {%
{\bf p}}_y=\hbar {\cal R}({\bf p}_y)f,
\end{equation}
where $f\in \Phi (E_q^2)$ and $\hbar $ is the Planck constant. The
quasi--regular representation (\ref{eq: b6}) is unitary with respect to the
scalar product (\ref{eq: b5}). Unitarity implies that the operators $\hat {%
{\bf p}}_x$ and $\hat {{\bf p}}_y$ are at least symmetric in the Hilbert
space $H(E_q^2)$. We assume that they have self--adjoint extensions. Using (%
\ref{eq: b12}) we have the following realization of the self--adjoint
operators $\hat {{\bf p}}_x$ and $\hat {{\bf p}}_y$ 
\begin{eqnarray}\label{eq: II13} 
\hat{{\bf p}}_x f({\bf z}_+, {\bf z}_-) =\frac{i\hbar}{\sqrt{2}} (
D_{q^{-2}}^{{\bf z}_+} f({\bf z}_+, {\bf z}_-) + 
D_{q^2}^{{\bf z}_-} f(q^{-2}{\bf z}_+, {\bf z}_-) ) \nonumber \\
\hat{{\bf p}}_y f({\bf z}_+, {\bf z}_-) =\frac{\hbar}{\sqrt{2}} (
D_{q^{-2}}^{{\bf z}_+} f({\bf z}_+, {\bf z}_-) - 
D_{q^2}^{{\bf z}_-} f(q^{-2}{\bf z}_+, {\bf z}_-) ), 
\end{eqnarray}
where $f\in \Phi (E_q^2)$ is regular function of ${\bf z}_{+}$ and ${\bf z}%
_{-}$ with the ordering (\ref{eq: b15}). The commutation relations satisfied
by the elements of the ''quantum mechanical q-phase'' space are 
\begin{equation}
\label{eq: II14}[\hat {{\bf x}},\hat {{\bf y}}]=-i\tanh \Lambda (\hat {{\bf x%
}}^2+\hat {{\bf y}}^2),\ \ \ [\hat {{\bf p}}_x,\hat {{\bf p}}_y]=-i\tanh
\Lambda (\hat {{\bf p}}_x^2+\hat {{\bf p}}_y^2)
\end{equation}
\begin{equation}
\label{eq: II15}[\hat {{\bf p}}_x,\hat {{\bf x}}]=[\hat {{\bf p}}_y,\hat {%
{\bf y}}]=i\hbar e^{\frac{2\Lambda }\hbar \hat {{\bf l}}},\ \ \ [\hat {{\bf p%
}}_x,\hat {{\bf y}}]=[\hat {{\bf p}}_y,\hat {{\bf x}}]=0,
\end{equation}
where $\hat {{\bf l}}$ is a function of $\hat {{\bf x}}$, $\hat {{\bf y}}$, $%
\hat {{\bf p}}_x$, $\hat {{\bf p}}_y$ which acts in $\Phi (E_q^2)$ as 
\begin{equation}
\label{eq: II16}\hat {{\bf l}}f({\bf z}_{+},{\bf z}_{-})=i\hbar [{\bf z}%
_{+}\partial _{{\bf z}_{+}}f({\bf z}_{+},{\bf z}_{-})-{\bf z}_{-}\partial _{%
{\bf z}_{-}}f(q^{-2}{\bf z}_{+},{\bf z}_{-})].
\end{equation}
It satisfies the commutation relations 
\begin{equation}
\label{eq: II16a}[\hat {{\bf l}},\hat {{\bf y}}]=-i\hbar \hat {{\bf x}},\ \
\ [\hat {{\bf l}},\hat {{\bf x}}]=i\hbar \hat {{\bf y}},
\end{equation}
\begin{equation}
\label{eq: II16b}[\hat {{\bf l}},\hat {{\bf p}}_x]=i\hbar \hat {{\bf p}}_y,\
\ \ [\hat {{\bf l}},\hat {{\bf p}}_y]=-i\hbar \hat {{\bf p}}_x.
\end{equation}
From the second commutation relations of (\ref{eq: II14}) we see that for
higher energies the non--commutativity between the momentum operators is
bigger. Note that the above commutation relations (\ref{eq: II16a}) and (\ref
{eq: II16b}) are the familiar quantum mechanical relations if $\hat {{\bf l}}
$ is interpreted as the angular momentum operator $\hat {{\bf p}}_x\hat {%
{\bf y}}-\hat {{\bf p}}_y\hat {{\bf x}}$ in $q\rightarrow 1$ limit.
Inspection of the commutation relations (\ref{eq: II15}) reveals that only
the section of the Hilbert space on which the expectation value of $\hat {%
{\bf l}}$ is non--negative is acceptable. Otherwise we may have smaller
(even zero ) than the usual non--commutativity between the momenta and
position operators which is unacceptable. This problem however is specific
to the two dimensional nature of the space. In three dimensions we would not
have such a unphysical situation.

Using (\ref{eq: II2}) and (\ref{eq: b15}) we can represent an arbitrary
element $f\in \Phi (E_q^2)$ as the infinite series 
\begin{equation}
f({\bf r},{\bf o})=\sum_{n=0}^\infty \sum_{j=-\infty
}^\infty c_{nj}{\bf r}^ne^{ij{\bf o}}
\end{equation}
subject to the square integrability condition (\ref{eq: b1}). The position $%
\hat {{\bf r}}$ and $\hat {{\bf o}}$ operators in the polar coordinate
system are given by 
\begin{equation}\label{eq: II17}
\hat {{\bf r}}f=f{\bf r},\ \ \ \hat {{\bf o}}f=f{\bf o},\ \
\ f\in \Phi (E_q^2).
\end{equation}
The radial momentum operator $\hat {{\bf p}}$ is expressed as 
\begin{equation}
\label{eq: II19}\hat {{\bf p}}=(\hat {{\bf p}}_{+}e^{i\hat {{\bf o}%
}}+e^{-i\hat {{\bf o}}}\hat {{\bf p}}_{-}).
\end{equation}
By the virtue of (\ref{eq: b12}) we get the following realization of $\hat {%
{\bf p}}$ 
\begin{equation}
\label{eq: II21}\hat {{\bf p}}f({\bf r},{\bf o})=\frac{i\hbar }{1+q}(qD_q^{%
{\bf r}}+q^{-1}D_{q^{-1}}^{{\bf r}})f({\bf r},{\bf o})+\frac{i\hbar }{(1+q)%
{\bf r}}f({\bf r},{\bf o})
\end{equation}
in the Hilbert space $H(E_q^2)$. Putting 
\begin{equation}
\label{eq: II22}f({\bf r},{\bf o})={\bf r}^{-1/2}\Psi ({\bf r},{\bf o})
\end{equation}
in (\ref{eq: II21}) we get the following realization of the radial momentum
operator $\hat {{\bf p}}$ 
\begin{equation}
\label{eq: II23}\hat p\Psi ({\bf r},{\bf o})=i\hbar \hat D_q^{{\bf r}}\Psi (%
{\bf r},{\bf o})
\end{equation}
in the Hilbert space $H^{\prime }(E_q^2)$. The scalar product in $H^{\prime
}(E_q^2)$ is 
\begin{equation}
\label{eq: II24}(\Psi ,\Psi ^{\prime })=\psi ({\bf r}^{-1}\Psi \Psi ^{\prime
*}),
\end{equation}
where $\psi (\cdot )$ is the invariant integral on $E_q(2)$ given by (\ref
{eq: b5}). In (\ref{eq: II23}) $\hat D_q^{{\bf r}}$ is q--derivative defined
as 
\begin{equation}
\label{eq: II25}\hat D_q^{{\bf r}}\Psi ({\bf r})=\frac{\Psi (q{\bf r})-\Psi
(q^{-1}{\bf r})}{(q-q^{-1}){\bf r}}.
\end{equation}
For the radial momentum operator we have the following commutation relations 
\begin{equation}
\label{eq: II26}[\hat {{\bf p}},\hat {{\bf r}}]=i\hbar e^{\frac \Lambda
\hbar \hat {{\bf l}}}\hat {{\bf a}},\ \ \ \ \ \ [\hat {{\bf p}},\hat {{\bf o}%
}]=0,
\end{equation}
where $\hat {{\bf a}}$ is the self--adjoint operator in $H^{\prime }(E_q^2)$
defined as 
\begin{equation}
\label{eq: II27}\hat {{\bf a}}\Psi ({\bf r},{\bf o})=\frac{q^{1/2}\Psi (q%
{\bf r},{\bf o})+q^{-1/2}\Psi (q^{-1}{\bf r},{\bf o})}{(q^{1/2}+q^{-1/2})}.
\end{equation}

\renewcommand{\theequation}{III.\arabic{equation}} \setcounter{equation}{0}

\section{Uncertainty Relations and Wave Packets}

Let $X_1$, $X_2$ be the self-- adjoint operators satisfying the commutation
relation 
\begin{equation}
[X_1, X_2] = iX_3. 
\end{equation}
The uncertainty relation for operators $X_1$, $X_2$ in the state $\psi$ is
given by 
\begin{equation}
\label{eq: III1}(\Delta X_1) (\Delta X_2) \geq \frac{1}{2} \mid (\psi ,\{
X_1^\prime, X_2^\prime \} \psi ) + i (\psi , X_3 \psi)\mid, 
\end{equation}
where 
\begin{equation}
\label{eq: III2}X_j^\prime = X_j - (\psi, X_j\psi ); \ \ \ \ \ j=1,2; 
\end{equation}
and $\{\cdot ,\cdot\}$ is the anticommutator. The equality sign in (\ref{eq:
III1}) holds if the condition 
\begin{equation}
\label{eq: III3}X_1^\prime \psi = c X_2^\prime \psi 
\end{equation}
is fulfilled. Here $c$ is some complex number. If $c$ is pure imaginary,
that is if $c^*=-c$ the first term of the right hand side of (\ref{eq: III1}%
) vanishes. We then have 
\begin{equation}
\label{eq: III4}(\Delta X_1) (\Delta X_2) \geq \frac{1}{2}\mid (\psi , X_3
\psi)\mid ; \ \ \ c^*=-c. 
\end{equation}
The corresponding state $\psi$ is called the wave packet.

\vspace{1cm}

(i) {\bf uncertainty relation for position operators}

Consider the position operators (\ref{eq: II17}) in the polar coordinates.
One can not have the regular solutions of (\ref{eq: III3}). It is not
difficult to check that the expectation value of the anticommutator $\{ \hat{%
{\bf r}}^\prime , \hat{{\bf o}}^\prime \}$ is zero on the vectors $f({\bf r}%
)\in \Phi (E_q^2)$ depending only of ${\bf r}$, that is 
\begin{equation}
\label{eq: III5}(\Delta \hat{{\bf r}}) (\Delta \hat{{\bf o}}) \geq \frac{%
\Lambda}{2} \mid (f , \hat{{\bf r}} f)\mid, 
\end{equation}
where the scalar product $(\cdot ,\cdot )$ is given by (\ref{eq: b5}). We
can minimize the right hand side of (\ref{eq: III5}) on the gaussian 
\begin{equation}
\label{eq: III6}f({\bf r})= \sqrt{\frac{2\log q^2}{q^2-1} } \frac{1}{%
\varepsilon} e^{-\frac{{\bf r}^2}{\varepsilon^2}} 
\end{equation}
as 
\begin{equation}
\label{eq: III7}(\Delta \hat{{\bf r}})(\Delta\hat{{\bf o}}) \geq \frac{%
\Lambda}{2} \sqrt{\frac{\pi}{2} } \varepsilon. 
\end{equation}
Making $\varepsilon$ very small we can decrease the uncertainty between the
position operators $\hat{{\bf r}}$ and $\hat{{\bf o}}$. Note that we cannot
put $\varepsilon =0$ since otherwise the wave packet (\ref{eq: III6})
becomes a non--regular function (precisely generalized function).

\vspace{1cm}

(ii) {\bf Uncertainty relation for momentum operators}

In this case the condition (\ref{eq: III3}) reads 
\begin{equation}
\label{eq: III8}(\hat{{\bf p}}_x -c\hat{{\bf p}}_y)f= (\langle {\bf p}%
_x\rangle -c \langle {\bf p}_y\rangle ) f. 
\end{equation}
For $c=i$ and $c=-i$ the above equation becomes 
\begin{equation}
\label{eq: III9}\hat{{\bf p}}_+ f_{\langle \hat{{\bf p}}_+\rangle}= \langle 
\hat{{\bf p}}_+\rangle f_{\langle \hat{{\bf p}}_+\rangle} 
\end{equation}
and 
\begin{equation}
\label{eq: III10}\hat{{\bf p}}_- f_{\langle \hat{{\bf p}}_-\rangle}= \langle 
\hat{{\bf p}}_-\rangle f_{\langle \hat{{\bf p}}_-\rangle} 
\end{equation}
respectively. By making use of (\ref{eq: b12}) we get the solutions 
\begin{equation}
\label{eq: III11}f_{\langle \hat{{\bf p}}_+\rangle}= N_+ e_{q^{-2}}^{ -i%
\frac{\langle \hat{{\bf p}}_+\rangle {\bf z}_+}{\hbar}} e_{q^{-2}}^{-i\frac{%
\langle \hat{{\bf p}}_-\rangle {\bf z}_-}{\hbar}}= N_+e_{q^{-2}}^{-i\frac{%
\langle \hat{{\bf p}}_+\rangle {\bf z}_+ + \langle \hat{{\bf p}}_-\rangle 
{\bf z}_-}{\hbar}}; \ \ c=i 
\end{equation}
and 
\begin{equation}
\label{eq: III12}f_{\langle \hat{{\bf p}}_-\rangle}= N_- e_{q^2}^{-i\frac{%
\langle \hat{{\bf p}}_-\rangle {\bf z}_-}{\hbar}} e_{q^{-2}}^{-i\frac{%
\langle \hat{{\bf p}}_+\rangle {\bf z}_+}{\hbar}}= N_-e_{q^2}^{-i\frac{%
\langle \hat{{\bf p}}_-\rangle {\bf z}_- + \langle \hat{{\bf p}}_+\rangle 
{\bf z}_+}{\hbar}}; \ \ c=-i. 
\end{equation}
Here $N_\pm$ are normalization constants. The above wave functions are the
q--analogues of the plane waves. They are not square integrable and
therefore do not belong to the Hilbert space $H(E_q^2)$. The wave packets
are given by 
\begin{equation}
f= \int db \ c(b) f_b, 
\end{equation}
where $c(b)$ is regular function of $b$ with the main value in the
neighborhood of ${\langle \hat{{\bf p}}_+\rangle}$ or ${\langle \hat{{\bf p}}%
_-\rangle}$ according to the choice (\ref{eq: III11}) and (\ref{eq: III12}).
We have the following uncertainty relations 
\begin{equation}
(\Delta \hat{{\bf p}}_x)(\Delta \hat{{\bf p}}_y ) \geq \frac{\tanh \Lambda}
{2} \mid (f , (\hat{{\bf p}}^2_x+\hat{{\bf p}}_y^2) f)\mid\sim \Lambda
\langle \hat{{\bf p}}_+ \hat{{\bf p}}_-\rangle. 
\end{equation}

\vspace{1cm}

(iii) {\bf $(\Delta\hat{ {\bf p} })(\Delta\hat{{\bf r}})$ uncertainty }

It is very difficult to solve (\ref{eq: III3}) for $\hat {{\bf p}}_x$, $\hat
{{\bf x}}$ and $\hat {{\bf p}}_y$, $\hat {{\bf y}}$ pairs. On the other hand
the commutation relations (\ref{eq: II14}), (\ref{eq: II15}) display a
circular symmetry. We then consider it is more meaningful to study the
uncertainty relation between the radial position $\hat {{\bf r}}$ and
momenta $\hat {{\bf p}}$ defined in (\ref{eq: II17}) and (\ref{eq: II19}).
For these operators the condition (\ref{eq: III3}) reads 
\begin{equation}
(\hat {{\bf p}}-c\hat {{\bf r}})\Psi ({\bf r},{\bf o})=(\langle \hat {{\bf p}%
}\rangle -c\langle \hat {{\bf r}}\rangle )\Psi ({\bf r},{\bf o})=d\Psi ({\bf %
r},{\bf o}),
\end{equation}
where $\Psi ({\bf r},{\bf o})$ is the vector from the Hilbert space $%
H^{\prime }(E_q^2)$ with scalar product (\ref{eq: II24}). By the virtue of (%
\ref{eq: II23}) we have 
\begin{equation}
\label{eq: iii}i\hbar \hat D_q^{{\bf r}}\Psi ({\bf r},{\bf o})-c{\bf r}\Psi (%
{\bf r},q^{-1}{\bf o})=d\Psi ({\bf r},{\bf o}),
\end{equation}
where $\hat D_q^{{\bf r}}$ is the q--derivative defined in (\ref{eq: II25}).
For the sake of simplicity let us consider the case $d=0$. The square
integrable solution of the equation (\ref{eq: iii}) exists for $c=-i\mid
c\mid q^{-j}$ which is given by 
\begin{equation}
\label{eq: iii1}\Psi _j({\bf r},{\bf o})=NE_{q^2}(-\frac{\mid c\mid q^{-j}}{%
(1+q)\hbar }{\bf r}^2)e^{ij{\bf o}},\ \ \ \ \ j=0,\ 1,\ 2,\ .\ .\ .,
\end{equation}
where $N$ is the normalization constant and $E_{q^2}(x)$ is the
q--exponential defined as 
\begin{equation}
E_{q^2}(x)=\sum_{n=0}^\infty \frac{x^n}{[n]_{q^2}!}.
\end{equation}
The q--factorial $[n]_{q^2}!$ is 
\begin{equation}
[n]_{q^2}!=[1]_{q^2}[2]_{q^2}\cdot \ \cdot \ \cdot [n]_{q^2},\ \ \ \
[n]_{q^2}=\frac{q^{2n}-q^{-2n}}{q^2-q^{-2}}.
\end{equation}
Note that the solution (\ref{eq: iii1}) is the eigenfunction of the angular
momentum operator $\hat {{\bf l}}$ corresponding to the e--value $j\hbar $ : 
\begin{equation}
\hat {{\bf l}}\Psi _j({\bf r},{\bf o})=\hbar j\Psi _j({\bf r},{\bf o}).
\end{equation}
Note that as we already discussed in the paragraph following (II.22)
physically only the chiral states with positive $j$ values are acceptable.
The minimal uncertainty relation between the radial momentum and position is 
\begin{equation}
(\Delta \hat {{\bf p}})(\Delta \hat {{\bf r}})\geq \frac \hbar 2\mid (\Psi
_j,\hat {{\bf a}}e^{\frac \Lambda \hbar \hat {{\bf l}}}\Psi _j)\mid ,
\end{equation}
where $\hat {{\bf a}}$ is the self--adjoint operator defined in (\ref{eq:
II27}). One can verify that in $q\rightarrow 1$ limit the first contribution
from the operator $\hat {{\bf a}}$ to the usual Heisenberg uncertainty
relation is of order $\Lambda ^2$. Therefore we can put $\hat {{\bf a}}=1$
in the right hand side of the above inequality to obtain the correction of
the order $\Lambda $. We then have 
\begin{equation}
(\Delta \hat {{\bf p}})(\Delta \hat {{\bf r}})\geq \frac \hbar 2\mid
1+\Lambda j\mid .
\end{equation}

\renewcommand{\theequation}{IV.\arabic{equation}} \setcounter{equation}{0}

\section{Free q--Semi--Classical Motion over $E_q^2$}

The dynamical model we have on the q--configuration space $E_q^2$ employs
classical time $t$. We now define the " quantum mechanical q--Hamiltonian''
which commutes with $\hat{{\bf p}}_x$, $\hat{{\bf p}}_y$ and $\hat{{\bf l}}$
: 
\begin{equation}
\label{eq: IV1}[\hat{{\bf H}}, \hat{{\bf p}}_x] = [\hat{\hat{{\bf H}}}, \hat{%
{\bf p}}_y] = [\hat{{\bf H}}, \hat{{\bf l}}]=0. 
\end{equation}
It is easy to check that the following operator 
\begin{equation}
\label{eq: IV2}\hat{{\bf H}}= \frac{1}{2m} (\hat{{\bf p}}^2_x + \hat{{\bf p}}%
^2_y ) e^{-\frac{2\Lambda }{\hbar}\hat{{\bf l}}} 
\end{equation}
satisfies the desired property (\ref{eq: IV1}), where $m$ is the mass of the
particle moving on $E_q^2$. We further assume that the commutator of any
operator with $\hat{{\bf H}}$ gives the time evolution of that operator.
Thus the time evolution of the position operators $\hat{{\bf x}}$, $\hat{%
{\bf p}}$ are 
\begin{equation}
\label{eq: IV3}\frac{d}{dt}\hat{{\bf x}}= \frac{1}{2m} [ \frac{(q+q^{-1})^2}{%
2}\hat{{\bf p}}_x + \frac{q^{-2}-q^2}{2}\hat{{\bf p}}_y ] +i\hbar^{-1}\hat{%
{\bf H}} [ \frac{(q-q^{-1})^2}{2}\hat{{\bf x}} + \frac{q^{-2}-q^2}{2}\hat{%
{\bf y}} ] 
\end{equation}
\begin{equation}
\label{eq: IV4}\frac{d}{dt}\hat{{\bf y}}= \frac{1}{2m} [ \frac{(q+q^{-1})^2}{%
2}\hat{{\bf p}}_y - \frac{q^{-2}-q^2}{2}\hat{{\bf p}}_x ] +i\hbar^{-1}\hat{%
{\bf H}} [ \frac{(q-q^{-1})^2}{2}\hat{{\bf y}} - \frac{q^{-2}-q^2}{2}\hat{%
{\bf x}} ]. 
\end{equation}
q--Hamiltonian ${\bf H}$ corresponding to (\ref{eq: IV2}) is given by 
\begin{equation}
\label{eq: IV5}{\bf H}= \frac{1}{2m} ({\bf p}^2_x + {\bf p}^2_y ) e^{-\frac{%
2\Lambda }{\hbar} {\bf l}}. 
\end{equation}
${\bf H}$ is the element of the q--phase space $E_q^2\otimes e_q^2$ defined
in the subsection (ii) of Sec. II. It acts in the Hilbert space $H(S\times S)
$.

We define the q--semiclassical Hamiltonian $H$ as the expectation value of $%
{\bf H}$ in the Hilbert space $H(S\times S)$ 
\begin{equation}
\label{eq: IV6}H=(u\otimes v,{\bf H}u\otimes v),
\end{equation}
where $u$ is the state from $H(S)$ on which the uncertainty between ${\bf x}$
and ${\bf y}$ is minimized and $v$ is the state on which the uncertainty
between ${\bf p}_x$ and ${\bf p}_y$ is minimized. They are given by 
\begin{equation}
u(\phi )=Ne^{-\frac{(\theta -\phi )^2}{2\Lambda }+i\frac{\log
\frac r{r_0}}\Lambda (\phi -\theta )}
\end{equation}
and 
\begin{equation}
v(\phi )=Ne^{-\frac{(\theta ^{\prime }-\phi )^2}{2\Lambda }+i%
\frac{\log \frac p{p_0}}\Lambda (\phi -\theta ^{\prime })},
\end{equation}
respectively. On this states we have 
\begin{equation}
\label{eq: IV9}x=(u,{\bf x}u)=r\cos \theta ,\ \ \ y=(u,{\bf x}u)=r\sin
\theta 
\end{equation}
and 
\begin{equation}
\label{eq: IV10}p_x=(v,{\bf p}_xv)=p\cos \theta ^{\prime },\ \ \ p_y=(u,{\bf %
p}_yu)=p\sin \theta ^{\prime }.
\end{equation}
It may be suggestive to sketch the quantum mechanical motion governed by the
q--semiclassical Hamiltonian 
\begin{equation}
\label{eq: IV11}H=\frac 1{2m}(p_x^2+p_y^2)e^{-\frac{2\Lambda }\hbar l},
\end{equation}
where $p_x$, $p_y$ and $l=xp_y-yp_x$ are the classical variables. The
Schr\"odinger equation corresponding to the above q--semiclassical
Hamiltonian, when it is written in plane--polar coordinates is 
\begin{equation}
-\frac{\hbar ^2}{2m}(\frac 1r\partial _r(r\partial _r)+\frac
1{r^2}\partial _\theta ^2)e^{-i\Lambda \partial _\theta }\psi (r,\theta
,t)=\frac i\hbar \psi (r,\theta ,t)
\end{equation}
solved by 
\begin{equation}
\label{eq: IV7}\psi _j(r,\theta ,t)=e^{-i\hbar \varepsilon }\frac 1{\sqrt{%
2\pi }}e^{ij\theta }J_j(\frac{\sqrt{2m\varepsilon }}\hbar e^{-\Lambda j/2}r),
\end{equation}
where $\varepsilon =\frac 1{2m}(p_x^2+p_y^2)$. The above Schr\"odinger
equation and its solution exhibit the peculiarity which is due to the
non--commutative nature of plane where semi q--classical motion we are
considering is originated. For states with $j=0$ the motion is exactly same
as usual one with $q=1$. However for particles with angular momentum $j$ the
radial distance is shorter for the chiral states corresponding to $j>0$ 
\begin{equation}
\label{eq: IV8}r\rightarrow e^{-\Lambda j/2}r=q^{j/2}r,\ \ \ \ for\ \ q<1
\end{equation}

\renewcommand{\theequation}{V.\arabic{equation}} \setcounter{equation}{0}

\section{Discussions}

Non--commutativity of the coordinates of the quantum plane is stronger for
larger radial distances (see (II.5) and (II.6)). Since the value of $q$ is
supposed to be very close to 1, the non--commutativity of $\hat {{\bf x}}$, $%
\hat {{\bf y}}$ coordinate operators is negligible if the expectation value
of $(\hat {{\bf x}}^2+\hat {{\bf y}}^2)$ is small. Thus we can say that in
the vicinity of the point where the observer is located the plane is almost
classical. For every observer the non--commutativity is same for the points
of plane corresponding to the same expectation value of $(\hat {{\bf x}%
}^2+\hat {{\bf y}}^2)$. In other words quantum plane is isotropic.

Non--commutativity of $\hat{{\bf p}}_x$ and $\hat{{\bf p}}_y$ on the other
hand can be appreciable at extremely high energies. We can then say that
this aspect of q--quantization should also manifest itself at very small
distances especially of the order of Planck scale.

Inspecting (III.23) we conclude that non--commutativity of $\hat {{\bf p}}_x$
and $\hat {{\bf x}}$ (and similarly $\hat {{\bf p}}_y$ and $\hat {{\bf y}}$
) differs from the usual one only for very large values of angular momentum
quantum number j.

\begin{center}
{\Large {\bf Appendix}}\cite{4}
\end{center}

\renewcommand{\theequation}{A.\arabic{equation}} \setcounter{equation}{0}

{\large A. Quantum Group $E_q(2)$ and it's Dual $e_q(2)$}

The Hopf algebra generated by ${\bf z}_{\pm }$ and ${\bf n}^{\mp 1}$ with
relations 
\begin{equation}
\label{eq: a1}{\bf z}_{+}{\bf z}_{-}=q^{-2}{\bf z}_{-}{\bf z}_{+},\ \ \ {\bf %
z}_{\pm }{\bf n}=q^2{\bf n}{\bf z}_{\pm },\ \ \ 
\end{equation}
involutions 
\begin{equation}
\label{eq: a2}{\bf n}^{*}={\bf n}^{-1},\ \ \ \ {\bf z}_{\pm }^{*}={\bf z}%
_{\mp }
\end{equation}
and group operations 
\begin{equation}
\label{eq: a3}\Delta ({\bf z}_{\pm })={\bf z}_{\pm }\otimes 1+{\bf n}^{\pm
1}\otimes {\bf z}_{\pm },\ \ \ \Delta ({\bf n}^{\pm 1})={\bf n}^{\pm
1}\otimes {\bf n}^{\pm 1},
\end{equation}
\begin{equation}
\label{eq: a4}\varepsilon ({\bf n}^{\pm 1})=1,\ \ \ \varepsilon ({\bf z}%
_{\pm })=0,\ \ \ 
\end{equation}
\begin{equation}
\label{eq: a5}S({\bf n}^{\pm 1})={\bf n}^{\mp 1},\ \ \ S({\bf z}_{\pm })=-%
{\bf n}^{\mp 1}{\bf z}_{\pm }.
\end{equation}
is called the quantum Euclidean group $E_q(2)$. The homomorphism 
\begin{equation}
\phi _K({\bf z}_{\pm })=0,\ \ \ \phi _K({\bf n}^{\pm 1})=t^{\pm 1}
\end{equation}
defines the quantum subgroup $U(1)$ of $E_q(2)$. We have the decomposition 
\begin{equation}
\label{eq: a6}E_q(2)=\sum_{i,j\in Z}\oplus E[i,j]
\end{equation}
of $E_q(2)$ with respect to it's subgroup $U(1)$. Here $E[i,j]$ are the
subspaces of $E_q(2)$ defined as 
\begin{equation}
\label{eq: a7}E[i,j]=\{f\in E_q(2):L_K(f)=t^i\otimes f;\ \ R_K(f)=f\otimes
t^j\},
\end{equation}
where 
\begin{equation}
\label{eq: a8}L_K=(\phi _K\otimes id)\circ \Delta ,\ \ \ \ R_K=(id\otimes
\phi _K)\circ \Delta .
\end{equation}
The quantum algebra $e_q(2)$ is the Hopf algebra generated by ${\bf p}_{\pm }
$ and ${\bf k}^{\pm 1}$ satisfying the relations 
\begin{equation}
\label{eq: a9}{\bf p}_{+}{\bf p}_{-}=q^{-2}{\bf p}_{-}{\bf p}_{+},\ \ \ \ \
\ {\bf p}_{\pm }{\bf k}=q^{\pm 2}{\bf k}{\bf p}_{\pm },
\end{equation}
involutions 
\begin{equation}
\label{eq: a10}{\bf p}_{\pm }^{*}={\bf p}_{\mp },\ \ \ \ {\bf k}^{*}={\bf k}
\end{equation}
and co--algebra operations 
\begin{equation}
\label{eq: a11}\Delta ({\bf p}_{\pm })={\bf p}_{\pm }\otimes 1+{\bf k}%
\otimes {\bf p}_{\pm },\ \ \ \Delta ({\bf k}^{\pm })={\bf k}^{\pm }\otimes 
{\bf k}^{\pm },
\end{equation}
\begin{equation}
S({\bf p}_{\pm })=-k^{-1}{\bf p}_{\pm },\ \ \ \ S({\bf k}^{\pm })={\bf k}%
^{\mp 1},
\end{equation}
\begin{equation}
\varepsilon ({\bf p}_{\pm })=0,\ \ \ \ \varepsilon ({\bf k})=1.
\end{equation}
The quantum algebra $e_q(2)$ and the quantum group $E_q(2)$ are in duality.
The duality bracket is given by 
\begin{equation}
\label{eq: a12}\langle {\bf p}_{-}^{n^{\prime }}{\bf p}_{+}^{k^{\prime }}%
{\bf k}^{j^{\prime }}\mid {\bf z}_{-}^n{\bf z}_{+}^k{\bf n}^j\rangle =\frac{%
i^{m+k}(1-q^2)^n(1-q^{-2})^k}{(q^2;q^2)_n(q^{-2};q^{-2})_k}q^{-2j^{\prime
}j}\delta _{nn^{\prime }}\delta _{kk^{\prime }}
\end{equation}
with ${j^{\prime }}$, $j\in Z$ and $n$, $n^{\prime }$, $k$ and $k^{\prime }$
being positive integers. Quantum subalgebra $u(1)$ of $e_q(2)$ corresponding
to the quantum subgroup $U(1)$ of $E_q(2)$ is generated by ${\bf k}^{\pm 1}$.

\vspace{1cm} \renewcommand{\theequation}{B.\arabic{equation}} 
\setcounter{equation}{0}

{\large B. Representation of $E_q(2)$ and it's Dual $e_q(2)$ }

We have the following realization of commutation relations (\ref{eq: a1}) 
\begin{equation}
\label{eq: b0}{\bf z}_{+}=r_0e^{-i\Lambda \frac d{d\phi }-i\phi },\ \ \ {\bf %
z}_{-}=r_0e^{-i\Lambda \frac d{d\phi }+i\phi },\ \ \ {\bf n}=e^{-2i\phi }
\end{equation}
( with $q=e^{-\Lambda }$ ) in the Hilbert space $H(S)$ of the square
integrable functions on the circle $S^1$. In the above equation $r_0$ is a
constant of dimension of length. In the orthonormal basis $e_j=\frac 1{\sqrt{%
2\pi }}e^{ij\phi }$, $-\infty <j<\infty $ we have 
\begin{equation}
\label{eq: b00}{\bf z}_{+}e_j=r_0q^{-j+\frac 12}e_{j-1},\ \ \ {\bf z}%
_{-}e_j=r_0q^{-j-\frac 12}e_{j+1},\ \ \ {\bf n}^{\pm }e_j=e_{j\mp 2}.
\end{equation}
Denote by $\Phi (E_q(2))$ the space of regular functions of ${\bf z}_{\pm }$
and ${\bf n}^{\pm 1}$ with finite norm 
\begin{equation}
\label{eq: b1}\mid \mid f\mid \mid <\infty ,\ \ \ \ \ \ f\in \Phi (E_q(2)),
\end{equation}
where 
\begin{equation}
\label{eq: b2}\mid \mid f\mid \mid =\sqrt{\psi (ff^{*})}
\end{equation}
with $\psi $ being the invariant integral on $E_q(2)$ defined as 
\begin{equation}
\label{eq: b3}\psi (f)=(1-q^2)\sum_{j=-\infty }^\infty (e_j,{\bf \xi }%
fe_j),\ \ \ {\bf \xi }={\bf z}_{+}{\bf z}_{-}.
\end{equation}
If $f$ is only the function of ${\bf \xi }$ the above expression can be
rewritten by means of the $q$--integral as 
\begin{equation}
\label{eq: b4}\psi (f)=\int_0^\infty f(\xi )d_{q^2}\xi .
\end{equation}
The space $\Phi (E_q(2))$ can be equipped with the scalar product 
\begin{equation}
\label{eq: b5}(f^{\prime },f)=\psi (f^{\prime }f^{*}).
\end{equation}
In a similar manner we define the space $\Phi (E_q^2)$. It consists of the
regular functions of coordinates ${\bf z}_{\pm }$ satisfying the square
integrability condition (\ref{eq: b1}). Completing $\Phi (E_q^2)$ in the
norm (\ref{eq: b2}) we arrive at the Hilbert space $H(E_q^2)$.

The comultiplication 
\begin{equation}
\label{eq: b6}\Delta :\Phi (E_q^2)\rightarrow \Phi (E_q(2))\otimes \Phi
(E_q^2)
\end{equation}
defines a left quasi--regular representation of the quantum group $E_q(2)$
in $\Phi (E_q^2)$. Since the scalar product in $\Phi (E_q^2)$ is defined by
means of invariant integral this representation is unitary.

The right representation ${\cal R}$ of the quantum algebra $e_q(2)$
corresponding to the left quasi--regular representation (\ref{eq: b6}) of $%
E_q(2)$ is given by 
\begin{equation}
\label{eq: b7}{\cal R}(\phi )f=(\phi \otimes id)\circ \Delta (f),\ \ \ \ \
f\in \Phi (E_q^2),
\end{equation}
where $\phi \in e_q(2)$ and 
\begin{equation}
(\phi \otimes id)\circ \Delta (f)=\langle \phi \mid f_j\rangle f_j^{\prime
},\ \ \ \ \ \Delta (f)=f_j\otimes f_j^{\prime }.
\end{equation}
Here $\langle \cdot \mid \cdot \rangle $ is duality bracket given in (\ref
{eq: a12}). The linear space $\Phi (E_q^2)$ is common invariant dense domain
for the set of linear operators ${\cal R}(\phi )$, $\phi \in e_q(2)$. Since
the representation (\ref{eq: b6}) of $E_q(2)$ is unitary we have 
\begin{equation}
\label{eq: b8}({\cal R}(\phi )f,f^{\prime })=(f,{\cal R}(\phi ^{*})f^{\prime
}),\ \ \ \ f\in \Phi (E_q^2).
\end{equation}
Some of the relations satisfied by the right representation ${\cal R}$ are 
\begin{equation}
\label{eq: b9}{\cal R}(\phi \phi ^{\prime })={\cal R}(\phi ^{\prime }){\cal R%
}(\phi )
\end{equation}
and 
\begin{equation}
\label{eq: b10}{\cal R}({\bf p}_{\pm })(ff^{\prime })={\cal R}({\bf p}_{\pm
})ff^{\prime }+{\cal R}({\bf k})f{\cal R}({\bf p}_{\pm })f^{\prime },
\end{equation}
\begin{equation}
\label{eq: b11}{\cal R}({\bf k})(ff^{\prime })={\cal R}({\bf k})f{\cal R}(%
{\bf k})f^{\prime }.
\end{equation}
We further have 
\begin{equation}
\label{eq: b12}{\cal R}({\bf p}_{+})f({\bf z}_{+},{\bf z}_{-})=iD_{q^{-2}}^{%
{\bf z}_{+}}f({\bf z}_{+},{\bf z}_{-}),\ \ \ {\cal R}({\bf p}_{-})f({\bf z}%
_{+},{\bf z}_{-})=iD_{q^2}^{{\bf z}_{-}}f(q^{-2}{\bf z}_{+},{\bf z}_{-}),
\end{equation}
\begin{equation}
\label{eq: b13}{\cal R}({\bf k})f({\bf z}_{+},{\bf z}_{-})=f(q^{-2}{\bf z}%
_{+},q^2{\bf z}_{-})
\end{equation}
and 
\begin{equation}
\label{eq: b14}{\cal R}({\bf p}_{+})f({\bf \xi })=iD_{q^2}^{{\bf \xi }}f(%
{\bf \xi }){\bf z}_{-},\ \ \ {\cal R}({\bf p}_{-})f({\bf \xi }%
)=iq^{-2}D_{q^{-2}}^{{\bf \xi }}f({\bf \xi }){\bf z}_{+}.
\end{equation}
The q--derivative employed in the above equation is defined as 
\begin{equation}
D_q^xf(x)=\frac{f(x)-f(qx)}{(1-q)x}.
\end{equation}
Note that ordering convention we employ is such that ${\bf z}_{+}$'s are
always on the left of ${\bf z}_{-}$'s, that is 
\begin{equation}
\label{eq: b15}f({\bf z}_{+},{\bf z}_{-})=\sum_{n,m=0}^\infty c_{nm}{\bf z}%
_{+}^n{\bf z}_{-}^m,
\end{equation}
where $c_{nm}$ are arbitrary complex coefficients chosen to satisfy the
condition (\ref{eq: b1}).

\end{document}